# Bounds of the Mertens Functions

Darrell Cox[a], Sourangshu Ghosh[b], Eldar Sultanow[c]

[a]Grayson County College
Denison, TX 75020
USA
dilycox@cableone.net

[b]Department of Civil Engineering,
Indian Institute of Technology Kharagpur, West Bengal, India
sourangshug123@gmail.com

[c]Potsdam University,
14482 Potsdam,
Germany
Eldar.Sultanow@wi.uni-potsdam.de

## ABSTRACT

In this paper we derive new properties of Mertens function and discuss about a likely upper bound of the absolute value of the Mertens function $\sqrt{\log(x!)} > |M(x)|$ when $x > 1$. Using this likely bound we show that we have a sufficient condition to prove the Riemann Hypothesis.

## 1. Introduction

We define the Mobius Function $\mu(k)$. Depending on the factorization of n into prime factors the function can take various values in $\{-1, 0, 1\}$

- $\mu(n) = 1$ if n has even number of prime factors and it is also square-free(divisible by no perfect square other than 1)
- $\mu(n) = -1$ if n has odd number of prime factors and it is also square-free
- $\mu(n) = 0$ if n is divisible by a perfect square.

Mertens function is defined as $M(n) = \sum_{k=1}^{n} \mu(k)$ where $\mu(k)$ is the Mobius function. It can be restated as the difference in the number of square-free integers up to $x$ that have even number of prime factors and the number of square-free integers up to $x$ that have odd number of prime factors. The Mertens function rather grows very slowly since the Mobius function takes only the value $0, \pm 1$ in both the positive and negative directions and keeps oscillating in a chaotic manner. Mertens after verifying all the numbers up to 10,000 conjectured that the absolute value of $M(x)$ is always bounded by $\sqrt{x}$. This conjecture was later disproved by Odlyzko and te Riele[1]. This conjecture is replaced by a weaker one by Stieltjes[2] who conjectured that $M(x) = O(x^{\frac{1}{2}})$. Littlewood[3] proved that the Riemann hypothesis is equivalent to the statement that for every $\epsilon > 0$ the function $M(x)x^{-\frac{1}{2}-\epsilon}$ approaches zero as x → ∞. This proves that the Riemann Hypothesis is equivalent to conjecture that $M(x) = O(x^{\frac{1}{2}+\epsilon})$ which gives a rather very strong upper bound to the growth of $M(x)$. Although there exists no analytic formula, Titchmarsh[4] showed that if the Riemann Hypothesis is true and if there exist no multiple non-trivial Riemann zeta function zeros, then there must exist a sequence $T_k$ which satisfies $k \leq T_k \leq k + 1$ such that the following result holds:

$$M_0(x) = \lim_{k \to \infty} \sum_{\substack{\rho \\ |\gamma|<T_k}} \frac{x^\rho}{\rho \zeta'(\rho)} - 2 + \sum_{n=1}^{\infty} \frac{(-1)^{n+1}}{(2n)!\, n\zeta(2n+1)} \left(\frac{2\pi}{x}\right)^{2n}$$

Where $\zeta(z)$ is the Riemann zeta function, and $\rho$ are all the all nontrivial zeros of the Riemann zeta function and $M_0(x)$ is defined as $M_0(x) = M(x) - \frac{1}{2}\mu(x)$ if $x \in Z^+$, $M(x)$ Otherwise (Odlyzko and te Riele[1])

## 2. New Properties of Mertens Function

Lehman[5] proved that $\sum_{i=1}^{x} M(\lfloor x/i \rfloor) = 1$. In general, $\sum_{i=1}^{x} M(\lfloor x/(in) \rfloor) = 1, n = 1, 2, 3, \ldots, x$ (since $\lfloor \lfloor x/n \rfloor /i \rfloor = (\lfloor x/(in) \rfloor)$. Let $R'$ denote a square matrix where element (i, j) equals 1 if $j$ divides $i$ or 0 otherwise. (In a Redheffer matrix, element $(i, j)$ equals 1 if $i$ divides $j$ or if $j = 1$. Redheffer[6] proved that the determinant of the matrix equals the Mertens Function $M(x)$.) Let $T$ denote the matrix obtained from $R'$ by element-by-element multiplication of the columns by $M\left(\left\lfloor\frac{x}{1}\right\rfloor\right), M\left(\left\lfloor\frac{x}{2}\right\rfloor\right), M\left(\left\lfloor\frac{x}{3}\right\rfloor\right) \ldots M(\left\lfloor\frac{x}{x}\right\rfloor)$. For example, the T matrix for x = 12 is

$$T = \begin{bmatrix} -2 & 0 & 0 & 0 & 0 & 0 & 0 & 0 & 0 & 0 & 0 & 0 \\ -1 & -1 & 0 & 0 & 0 & 0 & 0 & 0 & 0 & 0 & 0 & 0 \\ -1 & 0 & -1 & 0 & 0 & 0 & 0 & 0 & 0 & 0 & 0 & 0 \\ -1 & -1 & 0 & -1 & 0 & 0 & 0 & 0 & 0 & 0 & 0 & 0 \\ 0 & 0 & 0 & 0 & 0 & 0 & 0 & 0 & 0 & 0 & 0 & 0 \\ 0 & 0 & 0 & 0 & 0 & 0 & 0 & 0 & 0 & 0 & 0 & 0 \\ 1 & 0 & 0 & 0 & 0 & 0 & 1 & 0 & 0 & 0 & 0 & 0 \\ 1 & 1 & 0 & 1 & 0 & 0 & 0 & 1 & 0 & 0 & 0 & 0 \\ 1 & 0 & 1 & 0 & 0 & 0 & 0 & 0 & 1 & 0 & 0 & 0 \\ 1 & 1 & 0 & 0 & 1 & 0 & 0 & 0 & 0 & 1 & 0 & 0 \\ 1 & 0 & 0 & 0 & 0 & 0 & 0 & 0 & 0 & 0 & 1 & 0 \\ 1 & 1 & 1 & 1 & 0 & 1 & 0 & 0 & 0 & 0 & 0 & 1 \end{bmatrix}$$

**Theorem (1)**: $\sum_{i=1}^{x} M(\lfloor x/i \rfloor)\, i = A(x)$

**Proof:** Let us now take $A(x) = \sum_{i=1}^{x} \varphi(i)$ where $\varphi$ is Euler's totient function. Let $U$ denote the matrix obtained from $T$ by element-by-element multiplication of the columns by $\varphi(j)$. The sum of the columns of $U$ then equals $A(x)$. Now since $i = \sum_{d|i} \varphi(d)$ we can write $\sum_{i=1}^{x} M(\lfloor x/i \rfloor)\, i$ (the sum of the rows of U) equals $A(x)$.

By the Schwarz inequality, $A(x)/\sqrt{x(x+1)(2x+1)/6}$ is a lower bound of $\sqrt{\sum_{i=1}^{x} M(\lfloor x/i \rfloor)^2}$. $A(x) = \sum_{i=1}^{x} \varphi(i)$ Is further simplified by Walfisz[7] and Hardy and Wright[8] as

$$A(x) = \sum_{i=1}^{x} \varphi(i) = \frac{1}{2}\sum_{k=1}^{x} \mu(k) \left\lfloor\frac{x}{k}\right\rfloor \left(1 + \left\lfloor\frac{x}{k}\right\rfloor\right) = \frac{3}{\pi^2}x^2 + O(n(\log n)^{2/3}(\log \log n)^{4/3})$$

$$\sqrt{\sum_{i=1}^{x} M(\lfloor x/i \rfloor)^2} > \frac{A(x)}{\sqrt{\frac{x(x+1)(2x+1)}{6}}} > \frac{\frac{3}{\pi^2}x^2}{\sqrt{\frac{x(x+1)(2x+1)}{6}}} = \frac{3\sqrt{3}}{\pi}\frac{x^2}{\sqrt{x(x+1)(x+1/2)}}$$

This can be further simplified to

$$\sqrt{\sum_{i=1}^{x} M(\lfloor x/i \rfloor)^2 / x} > \frac{3\sqrt{3}}{\pi}\frac{1}{\sqrt{(1+1/x)(1+1/2x)}}$$

Taking limit of infinity on both the sides, we get

$$\lim_{x \to \infty} \sqrt{\sum_{i=1}^{x} M(\lfloor x/i \rfloor)^2 / x} > \frac{3\sqrt{3}}{\pi}$$

This shows that $\sum_{i=1}^{x} M(\lfloor x/i \rfloor)^2$ at large values of x is greater than $\frac{27}{\pi^2} x^2$.

Let $\Lambda(i)$ denote the Mangoldt function ($\Lambda(i)$ equals $log(p)$ if $i = p^m$ for some prime $p$ and some $m \geq 1$ or 0 otherwise). Mertens[9] proved that $\sum_{i=1}^{x} M(\lfloor x/i \rfloor) \log i = \psi(x)$ where $\psi(x)$ denotes the second Chebyshev function ($\psi(x) = \sum_{i \leq x} \Lambda(i)$).

**Theorem (2)**: $\sum_{i=1}^{x} M\left(\left\lfloor \frac{x}{i} \right\rfloor\right) \log(i) \sigma_0(i)/2 = \log(x!)$

**Proof:** Let $\sigma_x(i)$ denote the sum of positive divisors function ($\sigma_x(i) = \sum_{d|i} d^x$). Replacing $\varphi(j)$ with $log(j)$ in the U matrix gives a similar result.

Let $\lambda(n)$ denote the Liouville function ($\lambda(1) = 1$ or if $= p_1^{a_1}..p_k^{a_k}, \lambda(n) = (-1)^{a_1+..+a_k}$). $\sum_{d|n} \lambda(d)$ Equals 1 if n is a perfect square or 0 otherwise. Let $L(x) = \sum_{n \leq x} \lambda(d)$. Let us also assume $H(x) = \sum_{n \leq x} \mu(n) \log(n)$. $H(x)/(x \log(x)) \to 0$ as $x \to \infty$ and $\lim_{x \to \infty} M(x)/x - H(x)/(x \log(x))) = 0$. The statement $\lim_{x \to \infty} M(x)/x = 0$ is equivalent to the prime number theorem. Also, $\Lambda(n) = - \sum_{d|n} \mu(d) \log(d)$. (Apostol[10]).

The generalization of the Euler's totient function is Jordan totient function. Let it be denoted as $J_k(x)$ which is defined as number of set of k positive integers which are all less than or equal to n that will form a co-prime set of $(k + 1)$ positive integers together with $n$. Let us define $(x) = \sum_{i=1}^{x} J_k(i)$. It is known that $\sum_{d|n} J_k(d) = n^k$. Then we get the following theorem.

**Theorem (3):** $\sum_{i=1}^{x} M\left(\left\lfloor \frac{x}{i} \right\rfloor\right) i^k = B(x)$

$B(x)$ is expanded by McCarthy[11] to be

$$B(x) = \sum_{i=1}^{x} J_k(i) = \frac{n^{r+1}}{(r+1)\zeta(r+1)} + O(n^r)$$

We therefore get $\sum_{i=1}^{x} M\left(\left\lfloor \frac{x}{i} \right\rfloor\right) i^k = B(x) = \frac{n^{k+1}}{(k+1)\zeta(k+1)} + O(n^k)$.

Likewise we can derive some other similar relationships using the $T$ matrix that are as listed below:

**Theorem (4):** $\sum_{i=1}^{x} M\left(\left\lfloor \frac{x}{i} \right\rfloor\right) \sigma_k(i) = \sum_{i=1}^{x} i^k$ for $k \in Z^+$

**Theorem (5):** $\sum_{i=1}^{x} M\left(\left\lfloor \frac{x}{i} \right\rfloor\right)$ where the summation is over those $i$ values that are perfect squares equals $L(x)$

**Theorem (6):** $\sum_{i=1}^{x} M\left(\left\lfloor \frac{x}{i} \right\rfloor\right) \Lambda(i) = -H(x)$

**Theorem (7):** $\sum_{i=1}^{x} M\left(\left\lfloor \frac{x}{i} \right\rfloor\right) 2^{\omega(n)} = \sum_{i=1}^{x} |\mu(i)| \sim \frac{6}{\pi^2} x^2 + O(\sqrt{n}) = No. of\ Square\ Free\ Integers$

**Theorem (8):** $\sum_{i=1}^{x} M\left(\left\lfloor \frac{x}{i} \right\rfloor\right) d(n^2) = \sum_{i=1}^{x} 2^{\omega(i)}$ where $d(x)$ is the sum of all the divisors of $x$

**Theorem (9):** $\sum_{i=1}^{x} M\left(\left\lfloor \frac{x}{i} \right\rfloor\right) d^2(n) = \sum_{i=1}^{x} d(i^2)$

**Theorem (10):** $\sum_{i=1}^{x} M\left(\left\lfloor\frac{x}{i}\right\rfloor\right)\left(\frac{i}{\varphi(i)}\right) = \sum_{i=1}^{x} \frac{\mu^2(i)}{\varphi(i)}$

Similarly many other relationships can be found between various arithmetic functions and the Mertens Functions.

## 3. A Likely Upper Bound of |M(x)|

The following conjecture is based on data collected for x ≤ 500, 000.

**Conjecture (1):** $log(x!) > \sum_{i=1}^{x} M(\lfloor x/i \rfloor)^2 > \psi(x)$ when $x > 7$

By Stirling's formula, $log(x!) = x\, log(x) - x + O(log(x))$, since $log(x)$ increases more slowly than any positive power of $log(x)$, this is a better upper bound of $\sum_{i=1}^{x} M(\lfloor x/i \rfloor)^2$ than $x^{1+\varepsilon}$ for any $\varepsilon > 0$. This likely bound can be used to prove the Riemann Hypothesis since $\sum_{i=1}^{x} M(\lfloor x/i \rfloor) > |M(x)|$ and therefore we can write $\sqrt{log(x!)} > |M(x)|$. Since the growth of $\sum_{i=1}^{x} M(\lfloor x/i \rfloor)^2$ is lesser than $x^{1+\varepsilon}$ for any $\varepsilon > 0$. We can say

$$M(x)x^{-\frac{1}{2}-\epsilon} \to 0 \ as \ x \to \infty.$$

Figure 1 for a plot of $(x!), \sum_{i=1}^{x} M(\lfloor x/i \rfloor)^2$, and $\psi(x)$ for $x = 1, 2, 3, \ldots, 1000$.

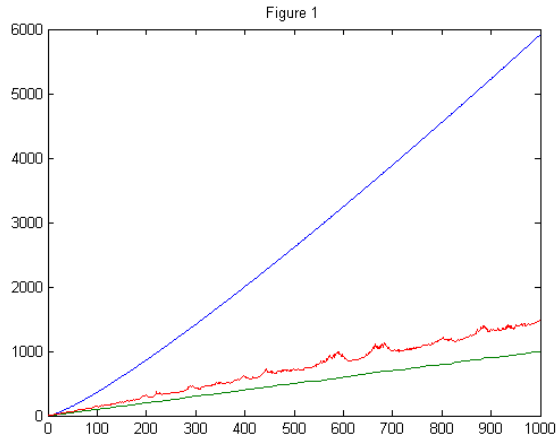

**Fig 1:** Plot of $log(x!), \sum_{i=1}^{x} M(\lfloor x/i \rfloor)^2$, and $\psi(x)$ for $x = 1, 2, 3, \ldots, 1000$

Let $j(x) = \sum_{i=1}^{x} M(\lfloor x/i \rfloor)^2$ where the summation is over $i$ values where $i|x$. Let $l_1, l_2, l_3$ denote the x values where $j(x)$ is a local maximum (that is, greater than all preceding $j(x)$ values) and let $m_1, m_2, m_3\ldots$ denote the values of the local maxima. The local maxima occur at x values that equal products of powers of small primes (Lagarias[12] discussed colossally abundant numbers and their relationship to the Riemann hypothesis). See Figure 2 for a plot of $l_i/(log(l_i)\, m_i)$, $m_i/l_i$, and $1/log(l_i)$ for $i = 1, 2, 3, \ldots, 772$ (corresponding to the local maxima for x ≤ 15, 000, 000, 000). ($M(x)$ Values for large x were computed using Del´eglise and Rivat's[13] algorithm.) The first two curves cross frequently, so there are $i$ values where mi is approximately equal to $l_i/\sqrt{log(l_i)}$.

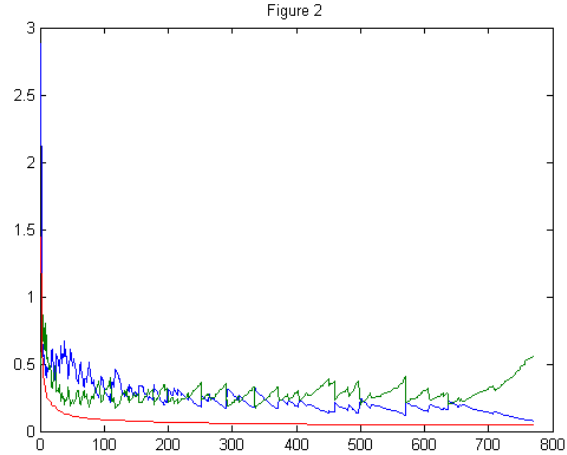

**Fig 2:** Plot of $l_i/(log(l_i)\, m_i)$, $m_i/l_i$, and $1/log(l_i)$ for $i = 1, 2, 3, \ldots, 772$

See Figure 3 for a plot of $j(x)$ and $\sum_{i=1}^{x} M(\lfloor x/i \rfloor)^2$ for x = 1, 2, 3, …, 10,000. See Figure 4 for a plot of $log(l_i), log(m_i), log(M(l_i)^2)$, and $log(m_i/\sigma_0(l_i))$ for i = 1, 2, 3, …, 772 (when $M(li) = 0, log(M(li)\, 2)$ is set to $-1$). See Figure 5 for a plot of $|M(l_i)|/\sqrt{l_i}$ for $i = 1, 2, 3, \ldots, 772$. The largest known value of $|M(x)|/\sqrt{x}$ (computed by Kotnik and van de Lune[14] for x ≤ 10$^{14}$) is 0.570591 (for $M(7,766,842,813) = 50,286$). The largest $|M(l_i)|/\sqrt{l_i}$ value for x ≤ 15,000,000,000 is 0.568887 (for $l_i = 7,766,892,000$). The largest known value of $|M(x)|/\sqrt{x}$ (computed by Kuznetsov[15] is 0.585767684 (for $M(11,609,864,264,058,592,345) = -1,995,900,927).)$

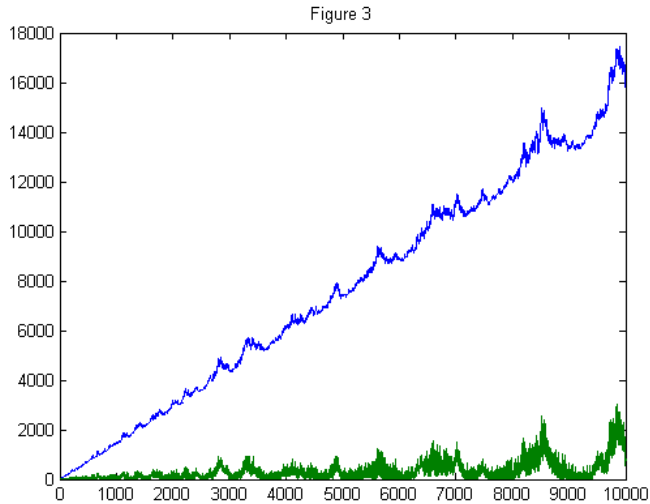

Fig 3: Plot of $j(x)$ and $\sum_{i=1}^{x} M(\lfloor x/i \rfloor)^2$ for x = 1, 2, 3, …, 10,000

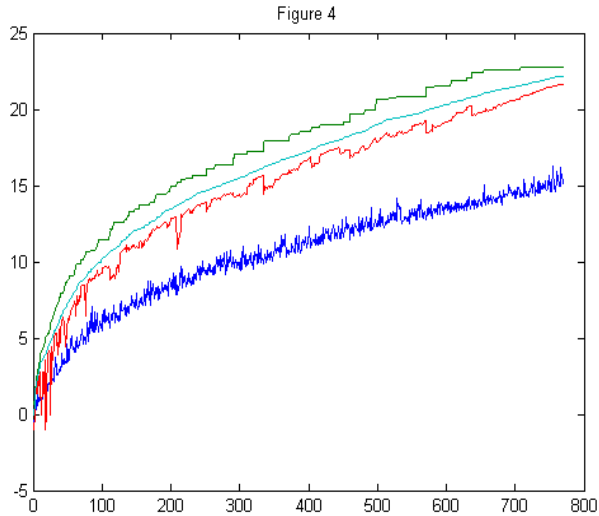

Fig 4: Plot of $log(l_i), log(m_i), log(M(l_i)^2)$, and $log(m_i/\sigma_0(l_i))$ for i = 1, 2, 3, ..., 772

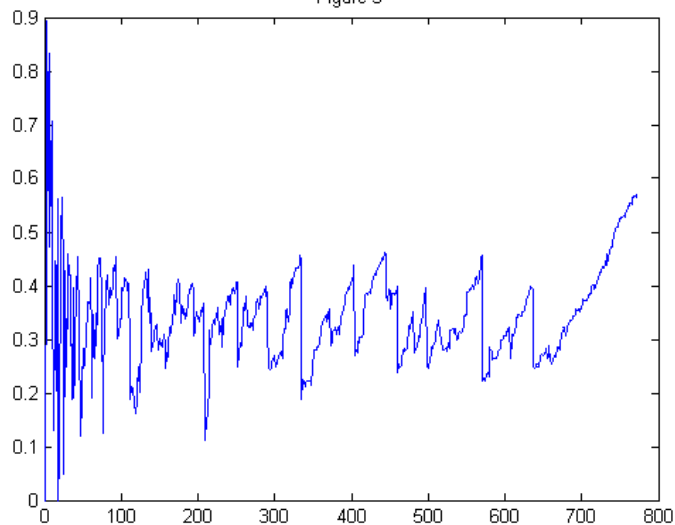

Fig 5: Plot of $|M(l_i)|/\sqrt{l_i}$ for $i = 1, 2, 3, \ldots, 772$

Let $l_i$ and $m_i$ be similarly defined for the function $\sigma_0(x)$. ($l_i$, $i = 1, 2, 3, \ldots$ are known as "highly composite" numbers. Ramanujan[16] initiated the study of such numbers. Robin[17] computed the first 5000 highly composite numbers.) Let $m_i'$ denote $j(l_i)$. See Figure 6 for a plot of $l_i/(log(l_i)m_i'), m_i'/l_i$, and $1/log(l_i)$ for $i = 2, 3, 4, \ldots, 160$ (corresponding to the local maxima for $x \leq 2,244,031,211,966,544,000$). (M(x) values for large x were computed using an algorithm similar to that used by Kuznetsov. The computations were done on an Intel i7-6700K CPU with 64 GB of RAM.) Although the first two curves cross frequently, $m_i'$ does not appear to converge to $l_i/\sqrt{log(l_i)}$.

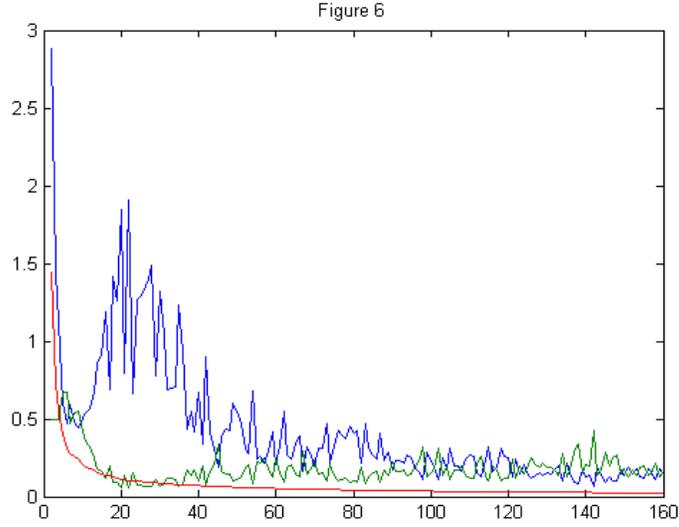

Fig 6: Plot of $l_i/(log(l_i)m_i')$, $m_i'/l_i$, and $1/log(l_i)$ for $i = 2, 3, 4, \ldots, 160$

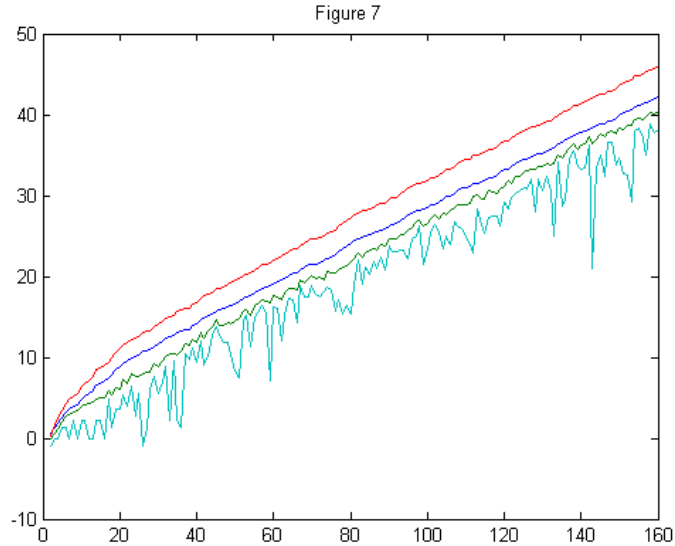

Fig 7: Plot of $log(l_i) + log(log(l_i))$, $log(l_i)$, $log(m_i')$, and $log(M(l_i)^2)$, for $i = 2, 3, 4, \ldots, 160$

See Figure 7 for a plot of $log(l_i) + log(log(l_i))$, $log(l_i)$, $log(m_i')$, and $log(M(l_i)^2)$, for $i = 2, 3, 4, \ldots, 160$ (when $M(l_i) = 0$, $log(M(l_i)^2)$, is set to $-1$). The vertical distance between the first and third curves appears to become roughly constant. See Figure 8 for a plot of $(log(l_i) + log(log(l_i)) - log(m_i')$, for i = 2, 3, 4... 160. See Figure 9 for a plot of $log(l_i) + \frac{1}{2}log(log(l_i))$, $log(\sum_{i=1}^{l_i} M(\lfloor l_i/i \rfloor)^2)$, and $log(l_i)$ for $i = 2, 3, 4, \ldots, 160$. $log(l_i) + \frac{1}{2}log(log(l_i))$ Is greater than $log(\sum_{n=1}^{l_i} M(\lfloor l_i/n \rfloor)^2)$ and $log(\sum_{n=1}^{l_i} M(\lfloor l_i/n \rfloor)^2)$ is greater than $log(l_i)$ for i > 4. This is evidence in support of Conjecture 1. See Figure 10 for a plot of $log(l_i) + \frac{1}{2}log(log(l_i)) - log(\sum_{n=1}^{l_i} M(\lfloor l_i/n \rfloor)^2)$ $for\ i = 2, 3, 4, \ldots, 160$.

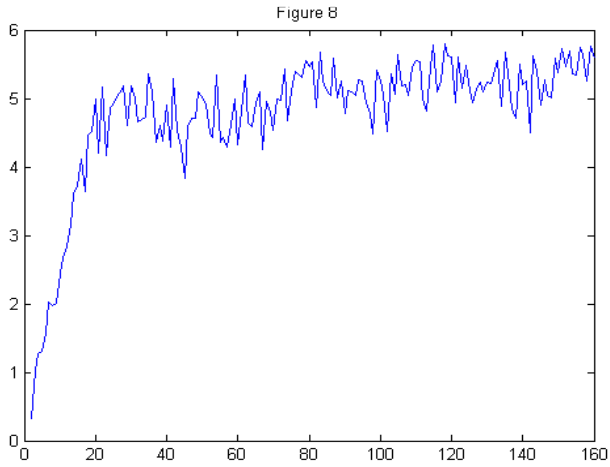

Fig: 8 for a plot of $(log(l_i) + log(log(l_i)) - log(m'_i)$, for i = 2, 3, 4... 160.

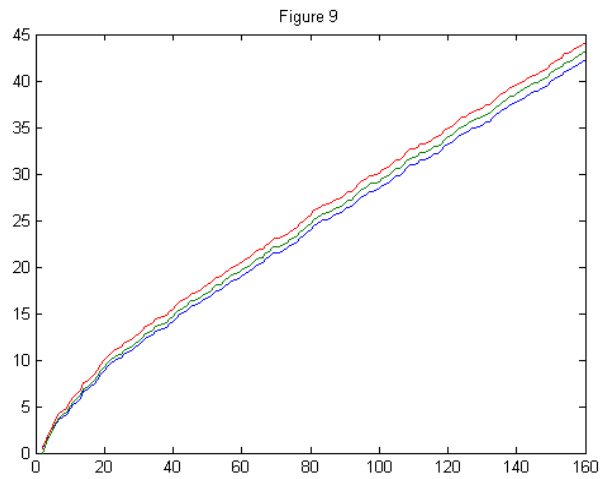

Fig 9: Plot of $log(l_i) + \frac{1}{2}log(log(l_i))$, $log(\sum_{i=1}^{l_i} M(\lfloor l_i/i \rfloor)^2)$, and $log(l_i)$ for $i = 2, 3, 4, \ldots, 160$

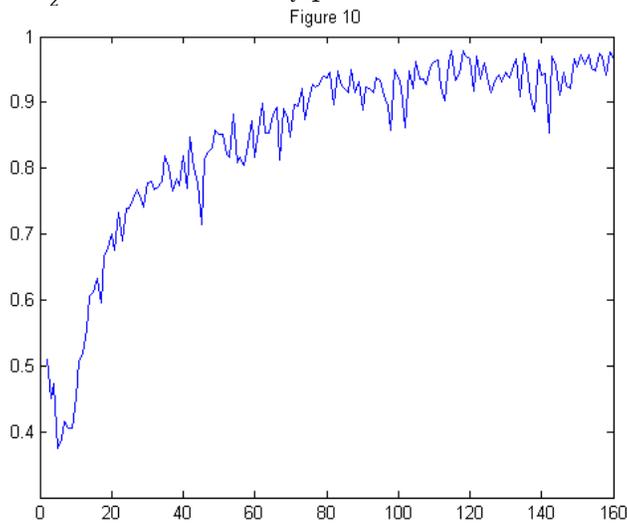

Fig 10: Plot of $log(l_i) + \frac{1}{2}log(log(l_i)) - log(\sum_{n=1}^{l_i} M(\lfloor l_i/n \rfloor)^2)$ for $i = 2, 3, 4, \ldots, 160$.

## 4. Conclusion

In this paper we derived new relations between Mertens function with a different arithmetic functions and also discussed about a likely upper bound of the absolute value of the Mertens function $\sqrt{\log(x!)} > |M(x)|$ when $x > 1$ with sufficient numerical evidence. More experimental evidence can be found in the work of Cox and Ghosh[18-19]. Using this likely upper bound we showed that we have a sufficient condition to prove the Riemann Hypothesis using the Littlewood condition.